\documentclass{ejm}

\usepackage[british]{babel}
\usepackage[utf8]{inputenc}
\usepackage[T1]{fontenc}
\usepackage{stmaryrd}

\usepackage{charter}

\usepackage{indentfirst}
\usepackage{xcolor}
\usepackage{verbatim}

\usepackage{amsmath}
\usepackage{amsfonts}	
\usepackage{amssymb	
           ,bbm
           ,units
           }
\usepackage{enumerate}

\usepackage{algorithm}
\usepackage[noend]{algpseudocode}
\usepackage{bm}

\makeatletter
\renewcommand{\ALG@beginalgorithmic}{\footnotesize}
\makeatother

\usepackage{hyperref}
\usepackage{wrapfig}
\usepackage{graphicx}
\usepackage{subfigure}

\usepackage[square,numbers,sort&compress]{natbib}
\numberwithin{equation}{section}
\numberwithin{figure}{section}
\numberwithin{algorithm}{section}
\numberwithin{table}{section}
 \usepackage[nodayofweek]{datetime}

\renewcommand*{\thefootnote}{\fnsymbol{footnote}}

\title[Hybrid PDE solver for data-driven problems and modern branching]{Hybrid PDE solver for data-driven problems and modern branching} 

\author[F. Bernal et al.]{Francisco Bernal$\,^{1}$, Gon\c calo dos Reis$\,^{2,3}$ \and  Greig Smith$\,^{2,4}$}
\affiliation{$^1\,$ CMAP - Centre de Math\'ematiques Appliqu\'ees,  Ecole Polytechnique, Route de Saclay, 91128 Palaiseau Cedex, FR.\\email: Francisco.Bernal@polytechnique.edu 
\\
$^2\,$ University of Edinburgh, School of Mathematics, Edinburgh, EH9 3FD, UK.\\ email: G.dosReis@ed.ac.uk
\\
$^3\,$ Centro de Matem\'atica e Aplica\c c$\tilde{\text{o}}$es (CMA), FCT, UNL, PT.
\\
$^4\,$ Maxwell Institute Graduate School in Analysis and its
	Applications (MIGSAA), University of Edinburgh, Edinburgh, EH9 3FD, UK.\\ email: G.Smith-13@sms.ed.ac.uk 
}

\renewcommand{\thefootnote}{\arabic{footnote}}

\newtheorem{theorem}{Theorem}[section]

\newtheorem{remark}[theorem]{Remark}

\newtheorem{assumption}[theorem]{Assumption}

\newcommand{\bE}{\mathbb{E}}

\newcommand{\bN}{\mathbb{N}}

\newcommand{\bR}{\mathbb{R}}

\newcommand{\bfx}{\mathbf{x}}

\newcommand{\cK}{\mathcal{K}}
\newcommand{\cL}{\mathcal{L}}

\newcommand{\ud}{\mathrm{d}}

\definecolor{darkgreen}{rgb}{0,0.35,0}

\newcommand{\1}{\mathbbm{1}}

\begin{document}

\selectlanguage{british}

\maketitle
\renewcommand*{\thefootnote}{\arabic{footnote}}

\begin{abstract}
The numerical solution of large-scale PDEs, such as those occurring in data-driven applications, unavoidably require powerful parallel computers and tailored parallel algorithms to make the best possible use of them. In fact, considerations about the parallelization and scalability of realistic problems are often critical enough to warrant acknowledgement in the modelling phase. The purpose of this paper is to spread awareness of the Probabilistic Domain Decomposition (PDD) method, a fresh approach to the parallelization of PDEs with excellent scalability properties. The idea exploits the stochastic representation of the PDE and its approximation via Monte Carlo in combination with deterministic high-performance PDE solvers. We describe the ingredients of PDD and its applicability in the scope of data science. In particular, we highlight recent advances in stochastic representations for nonlinear PDEs using branching diffusions, which have significantly broadened the scope of PDD. 

We envision this work as a dictionary giving large-scale PDE practitioners references on the very latest algorithms and techniques of a non-standard, yet highly parallelizable, methodology at the interface of deterministic and probabilistic numerical methods. We close this work with an invitation to the fully nonlinear case and open research questions.
\end{abstract}
{\bf Keywords:} Probabilistic Domain Decomposition, high-performance parallel computing, scalability, nonlinear PDEs, Marked branching diffusions, hybrid PDE solvers, Monte-Carlo methods.\bigskip


\noindent
{\bf 2010 AMS subject classifications:}\\
Primary: 65C05, 65C30, \quad 
Secondary: 65N55, 60H35, 91-XX, 35CXX


%
%
%

%

\section{Introduction}

Partial Differential Equations (PDEs) are ubiquitous in modelling, appearing in image analysis and processing, inverse problems, shape analysis and optimization, filtering, data assimilation and optimal control. They are used in Math-Biology to model population dynamics with competition or growth of tumours; or to model complex dynamics of movement of persons in crowds or to model (ir)rational decisions of players in games and they feature in many complex problems in Mathematical Finance. 
Underpinning all these applications is the necessity of solving numerically such equations either in bounded or unbounded domains. 

\emph{Deterministic Domain Decomposition.} The standard example of a Boundary Value Problem (BVP) is Laplace's equation with Dirichlet Boundary Conditions (BCs):
\begin{equation}
\label{F:Laplace}
\Delta u({\bf x})=0 \textrm{ if ${\bf x}\in\Omega\subset \bR^d$},\qquad
u({\bf x})= g({\bf x}) \textrm{ if ${\bf x}\in\partial\Omega$}.
\end{equation}

\begin{figure}[!ht] 
	\centering
\includegraphics[width=.375\columnwidth]{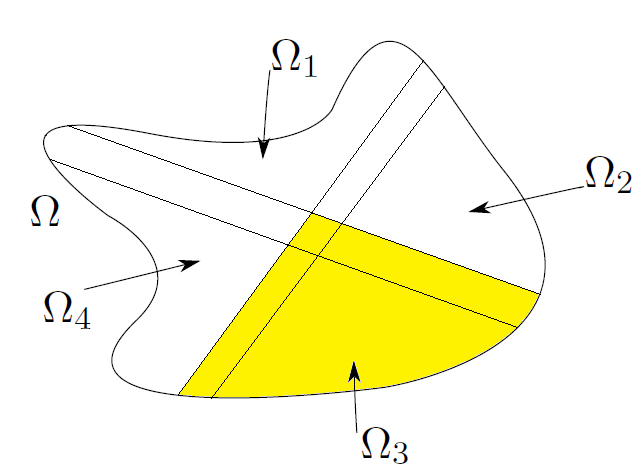}
\caption{Domain decomposition on an arbitrary domain $\Omega$, split into four overlapping subdomains $\Omega_{i}$ as required for Schwarz's alternating method. The subdomain $\Omega_3$ is highlighted.}
\label{fig: Basic Domain Decomposition}
\end{figure}

The large data sets involved in realistic applications nearly always imply that the discretization of a BVP such as (\ref{F:Laplace}) leads to algebraic systems of equations that can only be solved on a parallel computer with a large number (say $p>>1$) of processors. Not only does parallelization require multiple processors but also parallel algorithms. The classical Schwarz's alternating method was the first and remains the paradigm of such algorithms which we refer to as ``Deterministic Domain Decomposition'' (DDD) \cite{SmithBjorstadGropp2004}. While state-of-the-art DDD algorithms outperform Schwarz's alternating method in every respect, the latter nonetheless serves to illustrate the crucial difficulty they all face. The idea of Schwarz's algorithm is to divide $\Omega$ into a set of $p$ overlapping subdomains, and have processor $j=1,\cdots,p$ solve the restriction of the PDE to the subdomain, $\Omega_j$ see Fig.~\ref{fig: Basic Domain Decomposition}. 

Since the solution is not known in the first place, the BCs on the fictitious interfaces of $\Omega_j$ are also unknown, therefore, an initial guess has to be made in order to give processor $j$ a well-posed (yet incorrect) problem. The BCs along the fictitious interfaces of $\Omega_j$ are then updated from the solution of the surrounding subdomains in an iterative way until convergence. 

\emph{Inter-processor communication and the scalability limit of DDD.} Since the inter-processor communication involved in DDD's updating procedure is intrinsically sequential, it sets a limit to the scalability of the algorithm by virtue of the well known Amdahl's law. 

A simple illustrative example is as follows. A fully scalable algorithm would take half the time (say $T/2$) to run if the number of processors was doubled. If there is a fraction $\nu<1$ of the algorithm which is sequential, then the completion time will not drop below $T\nu$, regardless of how many processors are added. For instance, in Schwarz's method, if $5\%$ (i.e. $\nu=.05$) of the execution time of one given processor is lost by waiting for the artificial BCs to be ready, then the execution time could be shortened by at most a factor of 20 (with infinitely many processors; a factor of about 19 would already take over 1000 processors). In other words, Schwarz's alternating algorithm, or {\em any} DDD algorithm for that matter, cannot exploit the full capabilities of a parallel computer due to the idle time wasted in (essential) communication. 

We emphasize this point further by borrowing an example from David Keyes\footnote{See \href{http://www.mcs.anl.gov/research/projects/petsc-fun3d/Talks/bellTalk.ppt}{http://www.mcs.anl.gov/research/projects/petsc-fun3d/Talks/bellTalk.ppt}}. The Gordon Bell prizes are annually awarded to numerical schemes which achieve a breakthrough in performance when solving a realistic problem. In 1999, one such problem was the simulation of the compressible Navier-Stokes equations around the wing of an airplane. With $128$ processors, the winning code took $43$ minutes to solve the task. On the other hand, with $3072$ processors it took it $2.5$ minutes instead of $1.79$ as would have been the case with a fully parallelizable algorithm. The remaining $28\%$ of computer time were lost to interprocessor communication. At this point, adding more processors would have led to a faster loss of scalability. 

\emph{The Probabilistic Domain Decomposition (PDD) method.}
A conceptual breakthrough was achieved by Acebr\'on {\em et al.}  with the PDD method (or rather, the PDD framework) \cite{AcebronEtAl2005}, based on a previous, unpublished idea by Renato Spigler. PDD is the only domain decomposition method potentially free of communication, and thus potentially fully scalable. It does so by splitting the simulation in two separate stages, the first stage recasts the BVP into a stochastic formulation (via the so-called Feynman-Kac formula) which allows to compute the solution of the PDE at certain specific points in time/space. Thus, we can compute the ``true'' solution values of the DDD's fictitious interfaces for the $\Omega_j$'s. Therefore, the fictitious boundaries that were previously unknown are now known! 

Consequently, the subdomains are now completely independent of each other, the second stage then involves solving for the solution over the subdomains in a full parallel way. PDD calculations will be affected by two independent sources of numerical error: the subdomain solver and the statistical error of Monte Carlo simulations. 

PDD is currently well understood for the linear case, although recent advances in stochastic representations have opened the door for using PDD with nonlinear PDEs. A class of nonlinear PDEs strongly amenable to PDD are those whose solution can be represented by the so-called branching processes, as introduced by \cite{Watanabe1965}, \cite{Skorohod1964}, \cite{McKean1975} and recently extended by \cite{LabordereEtAl2014} and \cite{LabordereEtAl2016}. This methodology avoids backward regressions, the so-called ``Monte Carlo of Monte Carlo simulation'' problem, see Section \ref{Invitation} below or \cite[Section 3.1]{Labordere2012}. To illustrate the potential of branching in PDD,  numerical examples are worked out in Section \ref{Sec:Branching Diffusions KPP}. For general mildly nonlinear PDEs, the straightforward (and often efficient) approach of linearization and solution of each of the linear iterates, would be easily implementable with PDD, without resorting to nonlinear representations.

The general case of systems of nonlinear 2nd order parabolic/elliptic PDEs or fully nonlinear PDEs are, to the best of our knowledge, yet been addressed in the framework of PDD. These types of PDEs admit a probabilistic representation in terms of Forward Backward SDEs (FBSDEs \cite{ElKarouiPengQuenez1997}) and numerical methods for FBSDEs is currently the subject of extensive research.

\bigskip
The remainder of the paper is organized as follows.
\begin{itemize}
\item Section \ref{Overview} is an overview of PDD. 
We illustrate the connection between stochastic and BVPs in the linear case, and provide pointers for the numerical methods required. We stress the notion of balancing the various aspects of PDD in order to speed up the overall algorithm. The section closes with a survey of reported results on PDD.
\item In Section \ref{Sec:Branching Diffusions KPP}, nonlinear equations that can be represented probabilistically in terms of branching processes are discussed, and an account of recent developments is given. We give examples with convincing numerical tests using branching.
\item Section \ref{Invitation} is an invitation to new, more challenging problems than those tackled so far. In particular, the stochastic representation of systems of nonlinear parabolic equations with Forward Backward SDEs (FBSDEs) is succinctly discussed. We highlight open research problems.
\end{itemize}    

\section{An Overview of Probabilistic Domain Decomposition}\label{Overview}
\begin{figure}[!ht]
\centering
\includegraphics[width=0.90\columnwidth]{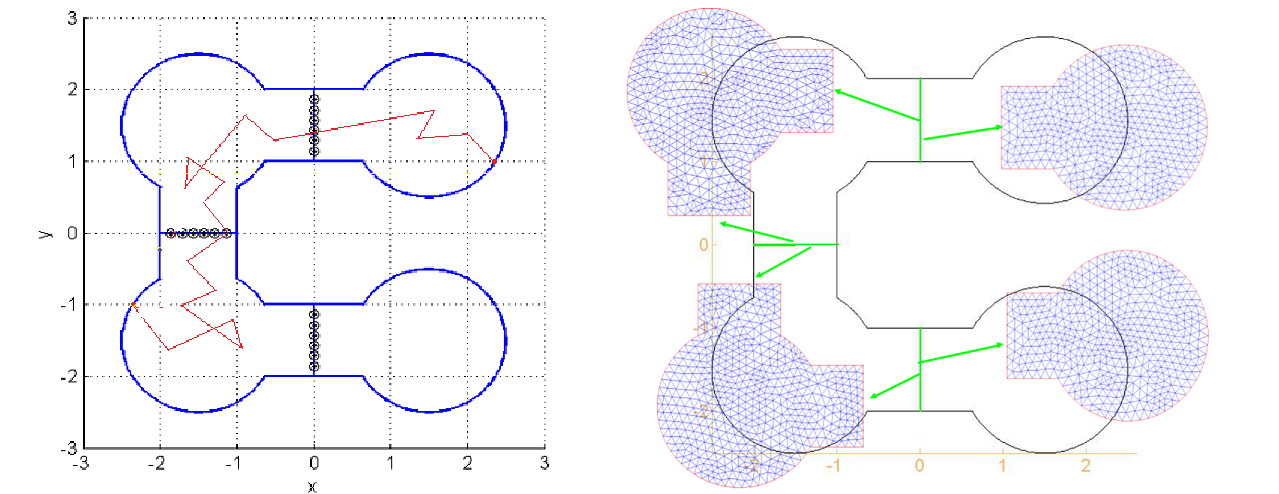}
\caption{An illustration of the two stages of PDD. (Left) First, the solution is computed on each of the interface nodes (black circles) between subdomains by numerically solving an appropriate SDE, which involves many independent realizations from the same interfacing node. Two such realizations are the trajectories shown in red. (Right) After the nodal solutions are known, the artificial interfaces (green segments) are provided with a Dirichlet BC by interpolation, and the BVPs on the various subdomains are decoupled. They can then be solved independently from one another with a deterministic method, such as FEM (meshes depicted).} 
\label{I:Figura1}
\end{figure}

An alternative to deterministic methods which is specifically designed to circumvent the scalability issue of DDD is the PDD method \cite{AcebronEtAl2005}, \cite{AcebronEtAl2005a}. First, the domain $\Omega$ under consideration is divided into non-overlapping subdomains. (Rather than overlapping ones, such as with Schwarz's alternating algorithm.) PDD then consists of two stages, firstly, the solution is calculated only on a set of interfacing nodes along the artificial interfaces, by solving the stochastic representation of the BVP with the Monte Carlo method, using the so-called Feynman-Kac representations. The stochastic representation is the crux of PDD and can be highly non-trivial. Nonetheless, it can also be extremely simple, for instance in (\ref{F:Laplace}), it is well known (see \cite{KaratzasShreve2012}) that the solution $u$ can be represented as
\begin{align}
\label{eq:ToyStochasticRepresentation} 
	u(\bfx) = \bE_{\bfx,0}\big[\, g({\bf X}_{\tau})\, \big]
\end{align}
where ${\bf X}$ is the solution to the simplest Stochastic Differential Equation (SDE)
\begin{align}
 \label{eq:SuperSimpleSDEforBVP} 
\text{for }t\geq 0\quad
	\ud {\bf X}_t & = \ud {\bf W}_t, \quad {\bf X}_0=\bfx\in\bR^d
	\quad\Leftrightarrow \quad {\bf X}_t = \bfx + {\bf W}_t
\end{align}
where $({\bf W}_t)_{t\geq 0}$ is a $d$-dimensional Brownian motion; ${\bf X}_{\tau}$ is the point on $\partial \Omega$ where the trajectory of ${\bf X}_t$ {\em first} hits the boundary and ${\bf X}$ starts at $t=0$ from $\bfx$; and lastly, $\bE[\cdot]$ is the usual expectation operator (in \eqref{eq:ToyStochasticRepresentation}, $\bE_{\bfx,0}[\cdot]$ emphasizes that the diffusion $({\bf X}_t)_{t\geq 0}$ starts at time $t=0$ from position $\bfx	$.) 
See Fig.\ref{I:Figura1} for an illustration.

The expected value above is approximated via a Monte Carlo method (the mean over many independent realizations of the SDE (\ref{eq:SuperSimpleSDEforBVP}), which can be carried out in a fully parallelizable way), introducing some statistical error. A numerical scheme is required to solve the SDE, we refer to it as ``stochastic numerics'' \cite{Kloeden_Platten}. 

By using the stochastic representation \eqref{eq:ToyStochasticRepresentation}, one can compute the equation's solution at each point at the artificial boundaries $\Omega_i \cap \Omega_j$.  It is then possible to reconstruct (approximately) the solution on the interfaces (interpolating the values at interfacing nodes with Chebyshev polynomials, say), so that the PDE restricted to each of the subdomains is now well posed. Now, the $p$ Laplace's equations on each of the $p$ subdomains are separate problems, and can be independently solved ``deterministically'', the second stage of PDD.  
Note that both stages in the PDD are \emph{embarrassingly parallel} by construction. We stress that this has been made possible only by the unique capability of stochastic representations of PDEs for solving a BVP at an isolated point (time-space) of the domain. On the other hand, DDD methods are matrix-based and deliver the solution everywhere; but by doing this are bound to move around information. 
\begin{remark}[Further features: Fault Tolerance and boundary design]
	Another feature of PDD is that it is, by construction, \emph{fault-tolerant}: the outage of an individual processor does not lead to substantial --- let alone fatal --- damage to the overall simulation. Given that modern parallel computers boast millions of processors, this is a definite plus (which stands in stark contrast to DDD algorithms).
	
	Lastly, it is up to the user to choose the artificial boundaries. These can be chosen in such a way as to lessen effects of kinks, sharp corners in the numerics. For instance, by better interplaying with the solutions being provided by the stochastic representations and Monte Carlo simulation. 
\end{remark}

%
%
%
%

\subsection{The ingredients of PDD}
PDD requires four ingredients: an interpolator, a subdomain solver, a stochastic representation of the PDE at hand, and an efficient numerical method for exploiting the latter with the Monte Carlo method. Let us briefly comment on them.

$\triangleright$ The {\em interpolator} scheme joins the (approximate) solutions at the interface nodes (computed with Monte Carlo) in order to yield an (approximate) 
Dirichlet BC on the whole interface. (We remark that, in time-dependent problems, the interface runs along the time component, as well as the spatial ones.) A flexible and accurate interpolation scheme well suited to PDD is the RBF interpolation --- see \cite{BernalAcebron2016a} and references therein. We note that less stringent approximation schemes than interpolation, such as ``denoising'' splines or least squares, are also possible.

$\triangleright$ The {\em subdomain solver} is any state-of-the-art numerical scheme to solve a PDE (or a BVP) on a finite domain, without parallelization. It takes the approximation to the solution on the boundary delivered by the interpolator as Dirichlet BCs and produces a solution everywhere within the subdomain. All the subdomain solutions are later ``glued together'' to yield a complete solution of the large-scale original problem. Common examples of solvers are finite differences, finite elements, and meshless methods. 

$\triangleright$ The {\em stochastic (or probabilistic)  representation of the PDE}, like that in \eqref{eq:ToyStochasticRepresentation}, is the pointwise characterization of the solution of a PDE (stationary or time-dependent) as the expected value of the functional of a system of stochastic differential equations (SDEs). The paradigm of stochastic representations is the well-known Feynman-Kac formula.  
More details will follow in Section \ref{Sec:Feynman Kac}.

$\triangleright$ {\em Efficient stochastic numerical methods} to exploit the stochastic representation. With the Monte Carlo approach, the expected value which yields the solution at an interfacing node is approximated as the average over many independent samples and this has to be repeated for every interfacing node. Contrary to deterministic methods, in the stochastic setting the errors are twofold, brought about by both discretization of time and by the substitution of the expected value by a sample mean. (One can say that the estimator of the expected value is biased due to the discretization of the time variable.) In order to attain a target accuracy $a$, both sources of errors must be balanced and this leads to lengthy simulations: the computational effort is ${\cal O}(a^{-4})$ if naive numerics are used (such as the Euler-Maruyama scheme plus boundary test for linear BVPs, see \cite{HighamMaoRojEtAl2013}). This cost can be substantially reduced (to around ${\cal O}(a^{-2})$ in some cases) with sophisticated, novel schemes. The new integrators for bounded stopped diffusions put forward by Gobet \cite{Gobet_HalfSpace} and Gobet and Menozzi \cite{GobetMenozzi2010} have doubled the convergence rate of the weak error associated to the interaction with the boundary, see also \cite{BernalAcebron2016}. For reflected diffusions, L\`epingle's method \cite{Gobet_HalfSpace} or Milstein's bounded methods \cite{Milstein_Tretyakov_Book} are recommended. Two further promising directions are the incorporation of Giles' Multilevel method \cite{HighamMaoRojEtAl2013}, \cite{GilesBernal2017} and the use of pathwise control variates for variance reduction \cite{BernalAcebron2016a}. In the later sections, pointers will be given to the numerics of nonlinear probabilistic representations. In any case, the numerics of SDEs (and generalizations thereof) are underdeveloped compared to those of deterministic equations. A positive aspect is that many new ideas being currently developed for SDEs in financial mathematics can be transplanted to PDD. 

Ideally, a numerical method tailored to the probabilistic representation of the given PDE exists, resulting in huge gains in efficiency. An example is the Walk on Spheres algorithm for the stopped Brownian motion, which is the probabilistic representation of Laplace's equation with Dirichlet BCs. This and similar, specific numerical schemes should be used whenever possible (see \cite{Chip} for references). 

\subsection{Probabilistic representation of linear BVPs}
\label{Sec:Feynman Kac}

We proceed now to give the probabilistic representation of linear BVPs (a generalization of the well known Feynman-Kac formula). Let us start by considering the general linear parabolic BVP of second order with mixed BCs:
\begin{equation}
\label{F:ParabolicBVP_mixedBCs}
\left\{
\begin{array}{ll}
\frac{\partial u}{\partial t}= \cL u + c({\bf x},t)u + f({\bf x},t) & \textrm{ if } t>0, {\bf x}\in\Omega,\\
u= p({\bf x}) & \textrm{ if } t=0, {\bf x}\in\Omega,\\
u= g({\bf x},t) & \textrm{ if } t>0, {\bf x}\in\partial\Omega_A,\\
\frac{\partial u}{\partial N}= {\varphi}({\bf x},t)u + \psi({\bf x},t) & \textrm{ if } t>0, {\bf x}\in\partial\Omega_R.
\end{array}
\right.
\end{equation}
The notation $\partial\Omega_A$ stands for a portion of the boundary where Dirichlet BCs are imposed, while on $\partial\Omega_R=\partial\Omega\backslash\partial\Omega_A$, BCs involving the normal derivative (Neumann or Robin) hold. In SDE literature $\partial\Omega_A$ and $\partial\Omega_R$ are known as absorbing and reflecting boundaries respectively. In (\ref{F:ParabolicBVP_mixedBCs}), $\cL$ is the 2nd order differential operator defined by 
\begin{align}
	\label{eq:DifferentialOperatorL}
	\cL:= \frac{1}{2} \sum_{i,j=1}^{d}A_{ij}(\bfx, t)\frac{\partial^2}{\partial x_{i} \partial x_{j}}+\sum_{i=1}^{d}b_{i}(\bfx, t)\frac{\partial}{\partial x_{i}},
\end{align}
where the matrix $A$ is positive definite and hence can be decomposed as $A:=[A_{ij}]=\sigma\sigma^T$ (for instance by Cholesky factorization). 
For future reference, we point out the importance of this operator $\cL$ in obtaining the (main) driving SDE, one can compare the coefficients of $\cL$ with $\bf X$ in \eqref{F:Milstein_sys} below.
It is assumed that $\varphi\leq 0$ and that the remaining coefficients and boundary in (\ref{F:ParabolicBVP_mixedBCs}) are smooth enough (in particular, the outward normal vector ${\bf N}$ to $\partial\Omega_R$ is well defined) so that an unique solution exists \cite{Constantini1998}, \cite{Freidlin1985}. 

In order to introduce the stochastic representation of (\ref{F:ParabolicBVP_mixedBCs}), let ${\bf X}_t,\,0\leq t\leq T$ be a diffusion (see \eqref{F:Milstein_sys} below) starting at ${\bf X}_0={\bf x}_0\in\Omega\subset{\mathbb R}^d,\,(d\geq 1)$, and consider the linear functional, in fact, a random variable,
\begin{align}
\label{F:Feynman-Kac}
\phi &= {\1}_{\{\tau\geq T\}} p({\bf X}_T)\Phi(T) + 
{\1}_{\{\tau<T\}}g({\bf X}_{\tau},T-\tau)\Phi(\tau)  \nonumber
\\
& \qquad \qquad
+
\int\limits_0\limits^{\min{(T,\tau})} f({\bf X}_t,T-t)\Phi(t)dt
+ \int\limits_0\limits^{\min{(T,\tau})} \psi({\bf X}_t,T-t)\Phi(t)d\xi_t, 
\end{align}
where
\begin{equation}
\label{F:Exponential}
\Phi(t)= \exp{\Big(\,\int_0^t c({\bf X}_s,s) ds + \int_0^t {\varphi}({\bf X}_s,s) d\xi_s\,\Big)}. 
\end{equation}
In (\ref{F:Feynman-Kac}) and (\ref{F:Exponential}),  
${\1}_{\{H\}}$ is the indicator function ($1$ if $H$ is true and $0$ otherwise); $\tau=\inf_{t}\{{\bf X}_t\in\partial\Omega_A\}$ is the ``first exit (or first passage) time'' from $\Omega$; which occurs at the ``first exit point'' ${\bf X}_{\tau}\in\partial\Omega_A$; $\xi_t$ is the ``local time'' (the amount of time the trajectory spends infinitesimally close to $\partial\Omega_R$)
; and all functions are assumed continuous and bounded. 
We want to calculate the expectation of $\phi$ in (\ref{F:Feynman-Kac}), which can also be expressed as
\begin{equation}
\label{F:Feynman-Kac-Milstein}
{\mathbb E}\big[\phi\big|{\bf X}_0={\bf x}_0\big]= {\mathbb E}\Big[\, q({\bf X}_{\tau})Y_{\tau}+Z_{\tau} \,\Big],
\ \textrm{s.~th. } \
q({\bf X}_{\tau})=\left\{
\begin{array}{ll}
g({\bf X}_{\tau},T-\tau), & \textrm{ if } \tau<T,\\
p({\bf X}_T), & \textrm{ if } \tau\geq T, 
\end{array}
\right.
\end{equation}
where the processes $({\bf X}_t,Y_t,Z_t,\xi_t)$ are governed by a set of stochastic differential equations (SDEs) driven by a $d-$dimensional Wiener process ${\bf W}_t$:  
\begin{equation}
\label{F:Milstein_sys}
\left\{
\begin{array}{ll}
d{\bf X}_t= {\bf b}({\bf X}_t,T-t)dt + \sigma({\bf X}_t,T-t)d{\bf W}_t -{\bf N}({\bf X}_t)d\xi_t & {\bf X}_0={\bf x}_0,\\
dY_t= c({\bf X}_t,T-t)Y_t dt + {\varphi}({\bf X}_t,T-t) Y_t d\xi_t & Y_0=1,\\ 
dZ_t= f({\bf X}_t,T-t)Y_t dt + \psi({\bf X}_t,T-t) Y_t d\xi_t  & Z_0= 0,\\
d\xi_t= {\1}_{\{{\bf X}_t\in\partial\Omega_R\}}dt & \xi_0=0.
\end{array}
\right.
\end{equation}
Then, formulas (\ref{F:Feynman-Kac}) and (\ref{F:Feynman-Kac-Milstein}) are the solution of the parabolic linear BVP, i.e.
\begin{equation}
u(t,{\bf x}_0)= {\mathbb E}\big[\phi\big|{\bf X}_0={\bf x}_0\big]= {\mathbb E}\Big[\, q({\bf X}_{\tau})Y_{\tau}+Z_{\tau} \,\Big]. 
\end{equation}
The representation in terms of the SDE system (\ref{F:Milstein_sys}), favoured by Milstein, is not the most usual one, but it can be programmed in a straightforward manner. If $c({\bf x})\leq 0$, elliptic BVPs can be formally derived from (\ref{F:ParabolicBVP_mixedBCs}), the Feynman-Kac formulas are then known as Dynkin's formulas. By taking $T\to \infty$, removing $p$, and dropping the time dependence from the surviving coefficients, the SDEs in (\ref{F:Milstein_sys}) are now autonomous as long as $\tau$ is finite. (This fails to happen, for instance, with a purely reflecting boundary $(\partial\Omega=\partial\Omega_R)$, when the solution $u({\bf x})$ is defined up to an arbitrary constant and requires one more compatibility condition \cite{Freidlin1985}.) In the case of purely Dirichlet BCs, $\xi_\cdot$ and the last equation in (\ref{F:Milstein_sys}) drop out. The representation of linear BVPs can be simplified in many cases, such as (\ref{eq:ToyStochasticRepresentation}) and (\ref{eq:SuperSimpleSDEforBVP}) for (\ref{F:Laplace}).


\subsection{An illustration of the effect of improved numerics}

Remarkable research into the numerical methods to solve representations of linear BVPs given above during the last fifteen years has meant that they can be solved today at a fraction of the cost required when PDD was introduced. To highlight the effect of improved stochastic numerics on the PDD methodology, let us take an example from \cite{BernalAcebron2016a}. There the authors study the BVP 
\begin{align}
\label{F:PDE_ejemplo}
\nabla^2 u + \frac{\cos{(x+y)}}{1.1+\sin{(x+y)}}\big(\frac{\partial u}{\partial x}+\frac{\partial u}{\partial y}\big)-\frac{x^2+y^2}{1.1+\sin{(x+y)}}u+f(x,y)= 0,
\end{align}
with $f(x,y)$ such that $u(x,y)=2\cos\big(2(y-2)x\big)+\sin\big(3(x-2)y\big)+3.1$ is the exact solution, as well as the Dirichlet boundary condition. The BVP domain can be seen in Figure \ref{I:Figura1}.

This problem was solved with the most current version of PDD (called IterPDD in \cite{BernalAcebron2016a}). IterPDD is a numerical suite which includes variance reduction techniques with iterative control variates based on nested, increasingly accurate global solutions. 

The speedup of a method A over a method B to solve a BVP is defined as
\begin{equation}
S(A,B)= \frac{\textrm{Time taken by method B}}{\textrm{Time taken by method A}}.
\end{equation}
Table \ref{Tabla_PDD} shows the speedup of IterPDD over the previous versions of PDD. Note that IterPDD retains all the advantages of any previous version of the PDD algorithm and hence it is not necessary to compare it against DDD solvers, since $S(IterPDD,DDD)=S(IterPDD,PDD)\times S(PDD,DDD)$. The theoretical speedup in Table \ref{Tabla_PDD} is the optimal speedup predicted by a sensitivity algorithm presented in \cite{BernalAcebron2016a}. Importantly, the latter is designed in such a way that it relies on fast warm up Monte Carlo sampling, and runs in parallel so as not to defeat the purpose of PDD. The agreement with the experimentally observed speedup is consistently conservative. 
\begin{table}
\centering
\begin{tabular}{lcc}
$a_0$ & \textrm{theoretical speedup} & \textrm{observed speedup (approx.)} \\
\hline
.04	  & 13.93  &  15 \\
.02 	& 28.34  &  29 \\
.01	  & 57.42  &  60 \\
.005	& 116.00 &  125 \\
.0025	& 233.84 &  250 \\
\hline
\end{tabular}
\caption{Acceleration of the most current version of PDD (IterPDD), with improved stochastic numerics, over the previous PDD algorithm. $a_0$ is the error tolerance on interfacing nodes. The BVP being solved is (\ref{F:PDE_ejemplo}).}
\label{Tabla_PDD}
\end{table}

The acceleration (speedup) grows larger as the target nodal accuracy $a_0$ (the largest admissible error of the numerical solution on the interface nodes, obtained via the probabilistic representation) decreases, approximately in an inversely proportional way. In summary, nearly every previously reported result using PDD, which were already faster than deterministic solvers, could be additionally accelerated by one to two orders of magnitude thanks to improved numerics. Further acceleration improvements are expected to follow in combination with Multilevel formulations \cite{GilesBernal2017}. 

\subsubsection*{Some applications of PDD}
 
The initial undertaking of the PDD programme was due to J. A. Acebr\'on and his research group \cite{AcebronEtAl2005}-\cite{AcebronRodriguez-Rozas2013}. Some realistic large-scale simulations carried out in supercomputers at the Barcelona Supercomputer Center and Rome's CASPUR proved the superiority of PDD over ScaLAPACK in terms of both total time and observed scalability when solving a nonlinear equations including KPP \cite{AcebronEtAl2010}, \cite{AcebronRodriguez-RozasSpigler2010}, \cite{AcebronRodriguezRozas2011}. 
Independently, Gobet and Mair\'e proposed a PDD-like domain decomposition scheme for the Poisson equation in \cite{GobetMaire2005}. Bihlo and Haynes have applied PDD to the generation of meshes for finite elements, a PDE-based task very well suited to stochastic representations \cite{BihloHaynes2016}, \cite{BihloHaynes2014}, \cite{BihloHaynesWalsh2015}.
Moreover, in \cite{AcebronRodriguez-Rozas2013}, the Vlasov-Poisson equations were tackled, which are of unquestionable interest in plasma physics. The PDD treatment of nonlinear PDEs like this one is further discussed in Section \ref{Sec:Branching Diffusions KPP}.



\section{Branching Diffusions and the KPP equation} 
\label{Sec:Branching Diffusions KPP}

The PDEs covered in Section \ref{Overview}, are all linear, unfortunately though PDEs arising in applications are typically not linear and standard SDE arguments are insufficient to address the general settings. That being said, one can derive stochastic representations for more general PDEs using so-called Forward Backward SDEs (FBSDEs), but FBSDEs are computationally expensive and difficult to handle (see Section \ref{Invitation}). However, \cite{RasulovRaimovaMascagni2010}, \cite{Labordere2012}, \cite{LabordereEtAl2014} and \cite{LabordereEtAl2016} have further developed the ``branching diffusions'' methodology as an efficient method to solve wider classes of nonlinear PDEs than the originally proposed by McKean. 

Branching diffusions allow one to tackle PDEs which are not linear without the FBSDE machinery. Although the classical results were somewhat restrictive on the form of the PDE (see \eqref{Eq:Classic branching diffusion} below), the recent developments allow one to address in an efficient way more general classes of PDEs. For example, when nonlinearities appear in the solution, or even gradients of the solution (see \eqref{Eq:SemiLinear with non-constant coefficients} \& Section \ref{sec:Agemarked} below). Throughout this section we supplement the theory with examples of PDEs that are amenable to branching and whose representations open the door for PDD as a viable algorithm for a far larger class of problems than previously available.

In Section \ref{Overview} we discussed stochastic representations for non-Dirichlet boundary conditions, unfortunately such boundary conditions have not been considered in this more general setting, but are future work for the authors. For the interested reader we mention recent analysis of the Monte Carlo Branching methodology to elliptic PDEs \cite{agarwal2017branching}.

\subsection{Reaction-Diffusion: KPP Equation} 

The branching diffusion idea is based on the works \cite{Skorohod1964}, \cite{Watanabe1965}, \cite{McKean1975}, to solve the KPP equation (in one spatial dimension with Dirichlet boundary condition),
\begin{align*}
\begin{cases}
\frac{\partial u}{\partial t}- \frac{\partial^2 u}{\partial x^2}-u(u-1)=0 \, , \quad x \in (x_{1},x_{2})\, , ~ x_{1},x_{2} \in \bR \, , ~ t>0, 
\\
u(x,0)=\psi(x) \, , \quad u(x_{1},t)=g(x_{1},t) \, , \quad u(x_{2},t) = g(x_{2},t) \, .
\end{cases}
\end{align*}
This was later generalized to allow for nonlinear PDEs such as,
\begin{align}
\label{Eq:Classic branching diffusion}
\frac{\partial u}{\partial t}-\cL u - c\left(\sum_{i=0}^{\infty}\alpha_{i}u^{i} -u \right)=0\, ,
\end{align}
where $c$ is a positive constant. Solutions to such nonlinear PDEs can be represented through so-called \emph{branching diffusion processes}. There has also been work for non integer powers of $u$, typically $u^{\alpha}$ for $\alpha \in [0,2]$, such process are referred to as \emph{super diffusions}, although we will not discuss these here to focus on more recent work, see \cite{Dynkin2004}, \cite{LabordereEtAl2014} and references therein for further details. Again the SDE process is governed by $\cL$, see \eqref{eq:DifferentialOperatorL}.

\begin{assumption}[Classical Branching Diffusions \cite{McKean1975}]
\label{assump:ClassicalKPP}
In the classical setting the stochastic representation of the above PDE relies on the following two conditions for the coefficients 
$\alpha$; $\alpha_{i} \ge 0$ and $\sum_{i=0}^{\infty}\alpha_{i}=1$.
\end{assumption}
To explain the stochastic representation to \eqref{Eq:Classic branching diffusion}, for point $(\bfx,t)$ we need to consider a set of particles. We start a particle at $\bfx$ at time $0$, this particle has a life time $\tau$ which is exponentially distributed with intensity $c$ (sometimes known as the ``branching rate''), we then simulate this particle until $\min(t, \tau,  \tau_{\partial \Omega})$, where $\tau_{\partial \Omega}$ is the particle hitting the spatial boundary, if $\tau$ is this smallest time, then the particle ``branches'' into $i$ particles with probability $\alpha_{i}$. These particles all start at the same point in space, however, they are all equipped with the own independent exponential random variable (``life time'') and driving Brownian motion. We then continue to simulate every particle until it has hit the space boundary or is alive at time $t$. We provide a detailed algorithm for branching diffusions (without space boundaries) in Section \ref{Sec:Reaction-Diffusion KPP}.

Consider the same initial and boundary value as above. The solution $u$ can then be written as the expected value of the product of the surviving particles or particles that hit the space boundary (see \cite{AcebronEtAl2010}),
\begin{align}
\label{Eq:Classic Branching with Boundary}
u(\bfx,t)= \bE_{\bfx,0}\left[ 
\prod_{i=1}^{N_{t}}\left(
\psi\big(X^{(i)}_{t}\big) \1_{\{t < \tau_{\partial \Omega}\}} 
+ g\big(X^{(i)}_{\tau_{\partial \Omega}},t-\tau_{\partial \Omega}\big) \1_{\{t \ge \tau_{\partial \Omega}\}} 
\right) 
\right] \, ,
\end{align}
where $X^{(i)}_{s}$ is the location of the $i$th particle surviving at time $s$ and $N_{t}$ is the number of surviving particles at time $t$ (a particle that hits the boundary is still ``alive'' at time $t$).

Although this set up allows for more general PDEs, it is not difficult to see that the computational complexity is greater for these types of problems (though still smaller than that for FBSDEs). Moreover, because of the nature of the branching, we are mainly confined to cases which are not ``overly'' nonlinear with time horizons that are not too long. Otherwise, we may end up with ``explosions'' whereby we are required to keep track of an extremely large number of particles, which will also tend to increase the number of branchings. We postpone discussion of this to Section \ref{Sec:Marked Branching Diffusions}.

\begin{figure}[h] 
	\label{fig:BranchingExample}
	\centering
\includegraphics[width=0.525\columnwidth,height=4cm]{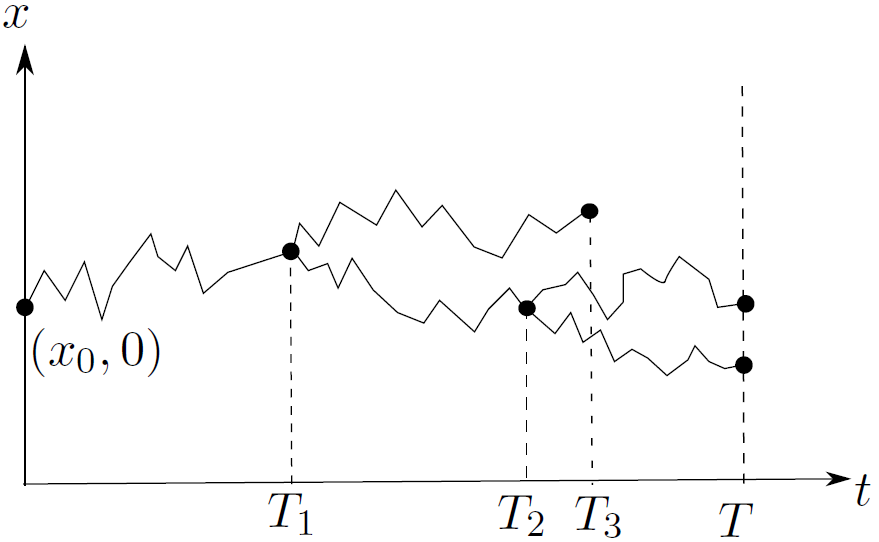}
	\caption{Representation of a branching diffusion, the first particle starts at point $x_{0}$ at time $t=0$, then after time $T_{1}$ branches into two particles. The first of these later dies at $T_{3}$ with no descendants, the other particle branches again in two particles at time $T_{2}$. Both of these (third generation) particles hits the terminal time, $T$ boundary.}
\end{figure}

To help give an idea of branching diffusions we illustrate a sample path for a given PDE. Let us consider the PDE (one spatial dimension),
\begin{align*}
\frac{\partial u}{\partial t}- \frac{\partial^2 u}{\partial x^2} - \frac{1}{2}\left(1+u^{2}\right) +u=0, \qquad u(x,0)=\psi(x),\quad x\in \bR \, , ~ t>0.
\end{align*}
From the work above there are two possible outcomes at a branching time (compare with \eqref{Eq:Classic branching diffusion}), either the particle dies, or splits into two descendants, the probability of each of these events is $\frac{1}{2}$. A possible trajectory (to compute $u(x_0,T)$) is shown in Figure \ref{fig:BranchingExample}.

\subsubsection{A PDD numerical example with the KPP Equation} 
\label{Sec:Reaction-Diffusion KPP}
Since we already have enough machinery to consider some interesting PDEs let us show an example which comes from the area of reaction diffusion equations, this class of equations are used in many areas such as physics and biology \cite{fisher1937wave}, \cite{FreiDosReis2013}, \cite{LionnetdosReisSzpruch2015}. 
For illustration purposes we consider the most well known of these, the Kolmogorov-Petrovsky-Piskunov (KPP) equation \cite{KolmogorovPetrovskyPiskunov1937}, linked to applications in plasma physics and ecology. The general form of the KPP is the following nonlinear parabolic PDE,
\begin{align*}
	\frac{\partial u}{\partial t}-D \frac{\partial^{2}u}{\partial x^{2}}+ru(1-u)=0 \, ,
	\qquad (x,t)\in \bR \times [0,\infty)
\end{align*}
where $D$ and $r$ are constants. Following \cite{AcebronEtAl2010}, taking $r=D=1$ and the initial value as,
\begin{align*}
	u(x,0)=\psi(x)=1-\Big(1+\exp\big(x/\sqrt{6}\big)\Big)^{-2}\, ,
\end{align*}
we can write the true solution of this problem as,
\begin{align*}
	u(x,t)=1-\Big(1+\exp\big(\dfrac{x}{\sqrt{6}}-\frac{5t}{{6}}\big)\Big)^{-2} \, .
\end{align*}
The advantage of such an example is that we can also compare the error of each method. Following \eqref{Eq:Classic Branching with Boundary}, the stochastic representation of the solution at $(x,t)$ is simply,
\begin{align*}
	u(x,t)
	= 
	\bE_{\bfx,0}\left[ \prod_{i=1}^{N_{t}}u\big(W^{(i)}_{t}+x,0\big) \right]
	=
	\bE_{\bfx,0}\left[ \prod_{i=1}^{N_{t}}\left(1- \left(1+\exp\left(\big(W^{(i)}_{t}+x\big)/\sqrt{6}\right)\right)^{-2}\right) \right] \, ,
\end{align*}
where $N_{t}$ is the number of particles at time $t$, and $W^{(i)}_{t}$ is the Brownian motion for particle $i$ at time $t$. Since the operator $\cL$ is $\partial^{2}/\partial x^{2}$, the process $X$ is Brownian motion started at point $x$ (compare with \eqref{F:ParabolicBVP_mixedBCs} and \eqref{F:Milstein_sys}).

\subsubsection*{Algorithm for the problem}
We solve this problem\footnote{MATLAB was used for the implementation. The simulations ran on a Dell PowerEdge R430 with four intel xeon E5-2680 processors. All polynomial interpolation were carried out using MATLAB's ``polyfit'' and the PDE is solved using MATLAB's ``pdepe'' function. To keep this consistent with the error we expect from Monte Carlo, we have set the relative and absolute error in the solver as $10^{-3}$.} over the domain $x \in [-2000,2000]$ for $t \in [0,1]$. For the PDD algorithm this corresponds to choosing points in space and solving them over a sequence of times to construct the artificial space time boundary, we create $p+1$ boundaries to obtain $p$ subdomains. Hence with five processors we split the domain in four subdomains $[-2000,-1000]$, $[-1000,0]$, $[0,1000]$ and $[1000,2000]$. This may not be the most optimal approach, but our goal is to show how even a simple PDD approach can dramatically improve the computational time required. We further calculate the solution at $11$ equally spaced time points (including $t=0$), which we denote by the set of $T_{i}$, satisfying $0=T_{0}<T_{1}< \dots<T=1$. Denote by $\Gamma$ the set of nodes at which we approximate the boundaries (true and artificial), hence for $D$ domain points, with $\Theta$ time points, $\Gamma$ is an $D$-by-$\Theta$ matrix. We denote by $\Gamma_{k}$ for $k \in \{1, \dots, D\}$ the column vectors of $\Gamma$, i.e. the approximation of the boundary at the $k$th domain point through time.
\begin{remark}[Pruning - Controlling the growth of the branching]
	A useful technique to control the branching is the ``pruning'', which is a method whereby we truncate the number of branches and comes from the observation that the pruned branches contribute very little to the expectation and are extremely computationally expensive. We do not discuss this further here, but direct the interested reader to \cite{AcebronEtAl2009}.
\end{remark}
It is clear from the PDE that all branching types will be two (the nonlinear part is a squared - compare with \eqref{Eq:Classic branching diffusion}). Moreover the driving process for the branching diffusions is a Brownian motion, this allows for us to take large step sizes since the simulation of Brownian motion is unbiased. Algorithm \ref{Alg:KPP} provides a summary of the method
\begin{algorithm}
	\label{Alg:KPP}
	\caption{Outline to solve the KPP equation using PDD over $p$ processors} 
	\begin{algorithmic}
		\State Split $\Omega$ into $(p-1)$ subdomains $\Omega_{i}$ and select set of nodes $\Gamma$ on the artificial boundary.
		\State Take a series points $T_{0}<T_{1}< \dots <T$ 
		and set the ``pruning'' level $P_{r}$.
		\For{Each $\Gamma_{k} \in \Gamma$}
		\For{$N$ Monte Carlo simulations}
		\While{$t <T$}
		\State Set $N_{p}$ as the number of particles and $t$ as the current time.
		\State Find $T_{i}$, such that $T_{i-1} \le t < T_{i}$ and simulate $t_{c} \sim \text{Exp}(c \times N_{p})$.
		\State Simulate each particle over $\min(t+t_{c},T_{i})$
		\If{we reach time $t+t_{c}$}
		\State Pick a particle (with equal probability) to branch into two
		\Else
		\State Evaluate $\psi^{(s)}(T_{i})=\prod_{j=1}^{N_{p}}\psi(X^{(j)}_{T_{i}})$.
		\EndIf{\textbf{end if}}
		\If{$N_{p}>P_{r}$}
		\State Restart this iteration of the for loop
		\EndIf{\textbf{end if}}
		\vspace{1mm}\EndWhile{\textbf{end while}}
		\vspace{1mm}\EndFor{\textbf{end for}}
		\vspace{1mm}\EndFor{\textbf{end for}}
		\For{Each $T_{i}$}
		\State Calculate $\Gamma_{k}(T_{i})=\frac{1}{N}\sum_{s=1}^{N}\psi^{(s)}(T_{i})$
		\EndFor{\textbf{end for}}
		\State Interpolate over each $\Gamma_{k}$ to create artificial boundaries and solve the PDE on each $\Omega_{i}$. 		
	\end{algorithmic}
\end{algorithm}

\subsubsection*{Error Analysis}

\begin{figure}[t]
	\centering
	\includegraphics[width=0.90\columnwidth,height=4.25cm]{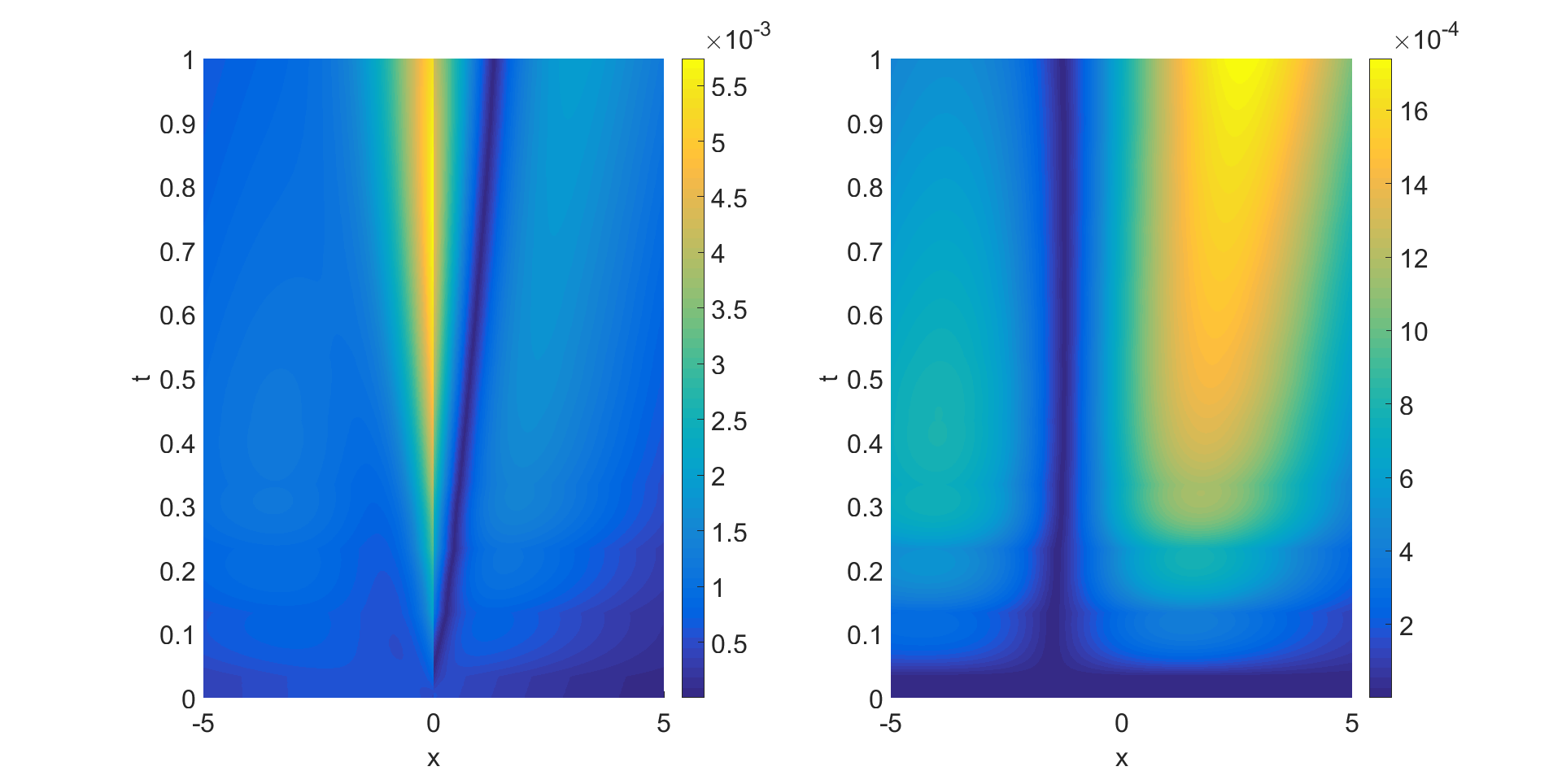}
	\vspace{-3mm}
	\caption{Absolute error in the region $[-5,5]$, for $N=10^5$ MC simulations (left), versus the pure PDE solver (right) with the PDE solver mesh set as $\Delta x=10^{-2}$ and $\Delta t=10^{-4}$.}
	\label{Fig:ErrorComparisonKPP}
\end{figure}

\begin{remark}[Straightforward implementation for Monte Carlo]
	\label{rem:cheapMCimplementation}
	For this example we used only a standard Monte Carlo implementation. Therefore, the results presented should be seen as the lower bound for what is achievable. It possible to improve the speed and accuracy of the algorithm by using methods touched on in Section \ref{Overview}. 
\end{remark}
Figure \ref{Fig:ErrorComparisonKPP} compares the absolute error of Monte Carlo versus using only the PDE solver over the domain $[-5,5]$. We observe the maximum error of the the pure PDE solver is approximately $1/3$ that of PDD. However, the largest error for PDD is concentrated in one region to the left hand side of the boundary and outside this small region the two errors are comparable.

\subsubsection*{Scalability and running times of the algorithm}
\begin{figure}[!hbt]
	\centering
	\subfigure
	{
		\includegraphics[width=0.51\columnwidth,height=4cm]{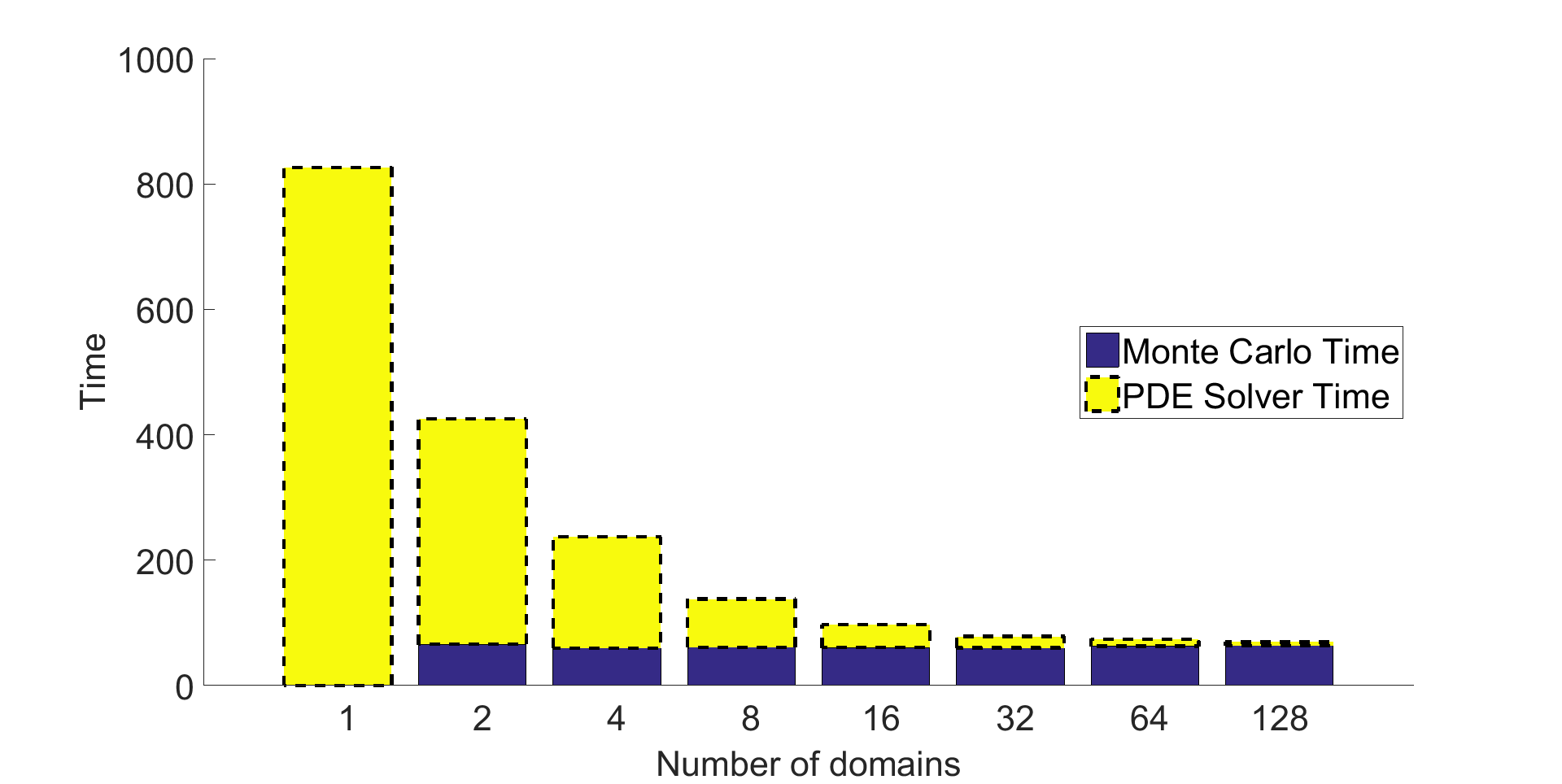}
	}
	\hspace{-1.05cm}
	\subfigure
	{
		\includegraphics[width=0.51\columnwidth,height=4cm]{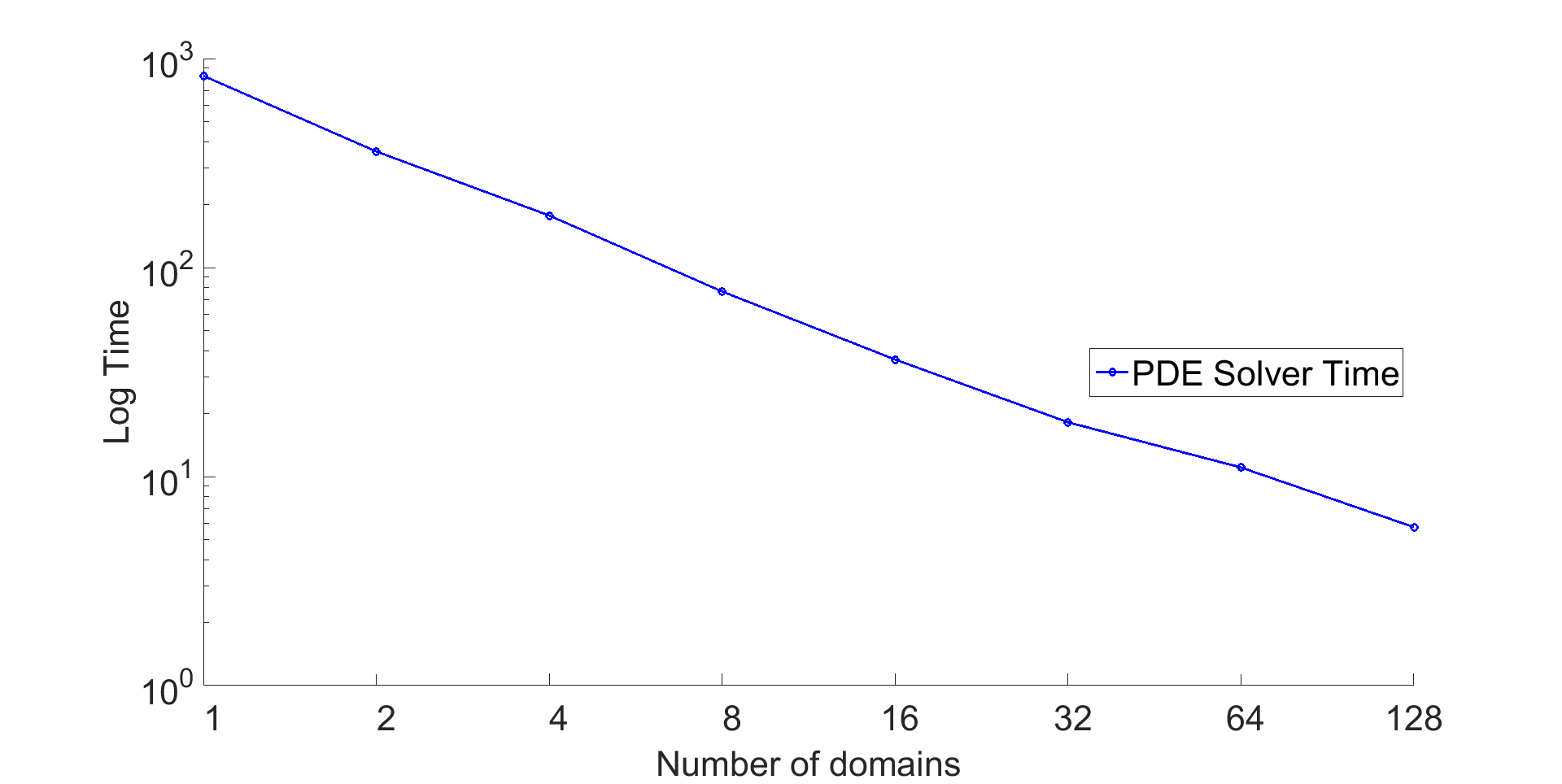}
	}
	\caption{On the left, Wall-clock running times (in seconds) for the PDE solver and PDD solver (Monte Carlo simulation and the PDE solver) as the number of subdomains (processors) increases. See Remark \ref{rem:cheapMCimplementation}. On the right, the running times in $\log$-scale showing perfect scalability.}
	\label{Fig:DomainDependentTimes}
\end{figure}

We conclude by considering how the time taken to solve the problem changes as more cores are used. We assume that we have access to the $p+1$ processors for $p$ subdomains. Since PDD consists of two main steps (Monte Carlo then PDE solver), we want to consider how these change as we increase the number cores. These results are presented in Figure \ref{Fig:DomainDependentTimes}. Note, to check the scaling we use the number of domains rather than the number of processors and PDD is not used when the number of domains is one. 

Figure \ref{Fig:DomainDependentTimes}, clearly shows the scaling capabilities of PDD. As one increases the number of processors the Monte Carlo time stays close to constant (as expected) and the PDE solver time drops dramatically as the number of processors increases. In fact, as is presented on the right graph, this decrease is constant with the number of subdomains with no apparent sign of levelling off; this follows from the no interprocessor communication.

\subsection{Marked Branching Diffusions}
\label{Sec:Marked Branching Diffusions}

The idea of \emph{marked branching diffusions} was developed by \cite{Labordere2012} as a solution to pricing credit valuation adjustments (CVAs) (see Section \ref{Sec:Pricing CVAs}), but since then has been generalized by \cite{LabordereEtAl2014} and \cite{LabordereEtAl2016}; similar ideas appeared in \cite{RasulovRaimovaMascagni2010}. Following \cite{LabordereEtAl2014}, we consider a terminal value PDE of the form (recall the time-inversion argument in Section \ref{Sec:Feynman Kac})
 \begin{align}
 \label{Eq:SemiLinear with non-constant coefficients}
 \dfrac{\partial u}{\partial t}+\cL u + c\left(\sum_{i=0}^{L}\alpha_{i}(\bfx,t)u^{i} -u \right)=0,
 \qquad 
 u(\bfx,T)=\psi(\bfx) \, .
 \end{align}
the notation follows that in Section \ref{Sec:Branching Diffusions KPP} with $L$ a positive integer.

We now discuss sufficient conditions for our stochastic representation to be the unique viscosity solution of this PDE. The technique in \cite{LabordereEtAl2014} is to use the fact that a large class of PDEs can be represented by FBSDEs, see \cite{PardouxPeng1992}. The idea is to derive sufficient conditions for the FBSDE to be the unique bounded viscosity solution to a PDE of the form \eqref{Eq:SemiLinear with non-constant coefficients} and use this to obtain the results on branching diffusions; here, FBSDEs are used only as a theoretical method to enable the branching argument. The key assumption is the following \cite[Assumption 2.2]{LabordereEtAl2014} denote the two functions for $s\in[0,\infty)$
\begin{align*}
l_{0}(s):= \sum_{k \ge 0} \|\alpha_{k}\|_{\infty} s^{k} 
\quad \text{and} \quad 
l(s):= c\Big(\, 
               \frac{l_{0}\big(s\|\psi\|_{\infty}\big)}{\|\psi\|_{\infty}}-s \,\Big),
\end{align*}
where $\| \cdot\|_{\infty}$ is the $L^{\infty}$-norm and the constants/functions follow from \eqref{Eq:SemiLinear with non-constant coefficients}, then,
\begin{assumption}[Assumptions for Marked branching]
	\label{Assump: Marked Branching Diffusion}
	$\phantom{123}${ }{ }
	\begin{enumerate}
		\item The power series $l_{0}$ has a radius of convergence $0<R \le \infty$. Moreover, the function $l$ satisfies one of the following conditions:
		\begin{enumerate}[(i)]
			\item $l(1) \le 0$,
			\item $l(1)>0$ and for some $\hat{s}>1$, $l(s)>0$, $\forall s \in [1,\hat{s})$ and $l(\hat s)=0$,
			\item $l(s)>0$ $\forall s \in [1, \infty)$ and $\int_{1}^{\bar{s}} \big(l(s)\big)^{-1}ds=T$, for some constant $\bar{s} \in (1,{R}/{\|\psi\|_{\infty}})$.
		\end{enumerate}
		\item The terminal function satisfies $\|\psi\|_{\infty}<R$.
	\end{enumerate}
\end{assumption}
Under the above assumption, the stochastic equations associated to \eqref{Eq:SemiLinear with non-constant coefficients} have a unique bounded solution and yield the unique viscosity solution for PDE \eqref{Eq:SemiLinear with non-constant coefficients} (see Proposition 2.3 and Theorem 2.13 of \cite{LabordereEtAl2014}). Thus the above assumption provides sufficient conditions for the branching diffusions to be used.

We are now most of the way to stating the marked branching diffusion representation. Of course, $\alpha$ needs not be a probability distribution any more, so let us introduce a probability distribution $q$, where $q$ is chosen so $q_{i}>0$ if $\alpha_{i} \neq 0$.
 It will beneficial for us to introduce notation to keep track of the particles. Let $T_{n}$ the $n$th branching time of the system where at time $T_{n}$ one of the particles branches into $k$ particles, with probability $q_{k}$. We denote by $I_{n}$ the number of descendants created at time $T_{n}$, hence $I_{n} \in \{0, \dots, L\}$. Denote by $M_{T-t}:=\sup\{n: t+T_{n} \le T\}$ the number of branches that occurred between time $t$ and $T$. Further denote by $\cK_{t}$ the index of all particles alive at time $t$, therefore $\cK_{T-t}$ corresponds to the index of all particles alive at time $T-t$. Finally, denote by $K_{n}$ the index of the particle that has branched at time $T_{n}$ (hence $K_{n} \in \cK_{T_{n}}$). The representation for the solution of the PDE \eqref{Eq:SemiLinear with non-constant coefficients} at $(\bfx,t)$ is given as
\begin{align*}
u(\bfx,t)=\bE_{\bfx,t}\left[ \prod_{k \in \cK_{T-t}}\psi({\bf X}^{(i)}_{T}) \prod_{n=1}^{M_{T-t}}\left(\frac{\alpha_{I_{n}}({\bf X}^{(K_{n})}_{t+T_{n}},t+T_{n})}{q_{I_{n}}}\right)\right] \, .
\end{align*}

\begin{remark}[Variance of the estimator]
	It is also of note that the specific choice of probability distribution $q$ (subject to the conditions above) does not change the representation, however, from the view point of variance reduction an optimal choice exists, see \cite{Labordere2012}.
	
	Although the set up here has been carried out using exponential random variables as the lifetime of particles, more recent work suggest that other distributions yield smaller variances, see \cite{LabordereEtAl2016}, \cite{DoumbiaEtAl2016}, \cite{Warin2017} and \cite{BouchardEtAl2016}.
\end{remark}

\subsubsection{More general PDEs through approximation}
\label{sec:obstaclePDE}
The stochastic representations shown in this section can be used to deal with many different types of PDEs. Namely, \emph{these representations allow us to consider PDEs whose source term $f$ can be well approximated by polynomials}. An example of this is considered in Section \ref{Sec:Pricing CVAs}. 

Although we mention CVAs, many other financial contracts rely on maximums between two entries (so called obstacle PDEs), these appear in American options for example. Computing the solution to certain obstacle PDEs in stochastics is deeply related to optimal stopping problems. We write such PDEs in their variational formulation: for a exercise payoff $\varphi$ the price function $u$ solves 
\begin{align*}
	\max\Big( \frac{\partial u}{\partial t} + \cL u,\varphi(\bfx) - u \Big) = 0, \qquad  u(\bfx,T) = \varphi(\bfx),
	\qquad (\bfx,t)\in \bR^d\times [0,T]	.
\end{align*}
In \cite{benth2003semilinear}, the authors showed how this PDE can be converted into a semilinear PDE 
\begin{align*}
	\frac{\partial u}{\partial t} + \cL u = \1_{\{ \varphi(\bfx)\geq u\}} \cL \varphi,\qquad u(\bfx,T)=\varphi(\bfx),
	\qquad (\bfx,t)\in \bR^d\times [0,T]	.
\end{align*}
A stochastic representation for these equations is well-understood and the methodology we have presented so far can be applied to them, see \cite[Section 2.3]{Labordere2012}. Moreover, \cite{BossyEtAl2015} use the polynomial representation of trigonometric functions to calculate the nonlinear Poisson-Boltzmann equation.

\subsubsection{Numerical example with Option Pricing: Credit Valuation Adjustment (CVA)} 
\label{Sec:Pricing CVAs}

The next example we consider is a different problem than the KPP above, but highlights the marked branching diffusion and its application scope very well. This example concerns option (derivative) pricing with Credit Valuation Adjustments (CVAs). Since the financial crisis in 2008, risk, especially default risk (the risk that a counter party fails to meet future obligations) has been at the forefront of many policy decisions and regulation changes for financial firms. CVAs play a crucial role here by adjusting the price of a derivative (a financial contract) with the knowledge that the counter party may default. \cite{GuyonLabordere2013} give a derivation of the PDEs arising in this problem. The specific form of the PDE depends on the way in which one chooses to model the mark-to-market value (the fair value) of the derivative at the time of default, but we obtain PDEs of the form,
\begin{align*}
	\dfrac{\partial u}{\partial t}+\cL u + r_{0}u+r_{1}u^{+}=0,
	\qquad 
	u(x,T)=\psi(x),\qquad x \in\bR
\end{align*}
where the $r_{i}$ are functions of $x$ and $t$ and we denote by 
$y^{+}:= \max(0,y)$.

\begin{figure}[!htb]
	\label{Fig:CVAPlot}
	\centering
	\includegraphics[width=0.75\columnwidth,height=4.5cm]{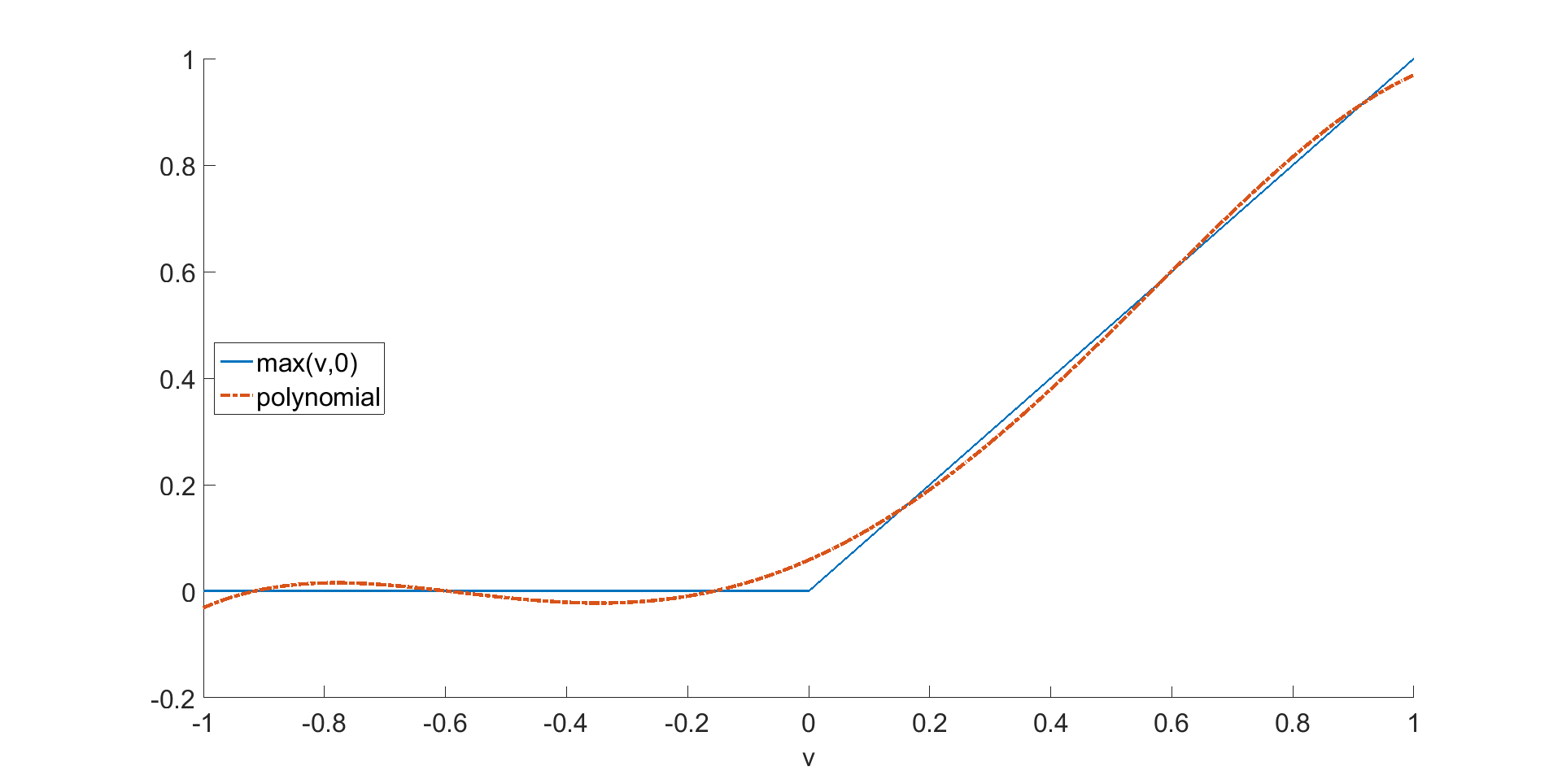}
	\caption{Approximation of $\max(v,0)$ over $v\in[-1,1]$ by a $4$th-order polynomial.}
\end{figure}

Although at first glance such a PDE seems beyond the scope of the stochastic representations in Section \ref{Sec:Marked Branching Diffusions}, \cite{Labordere2012} provided a simple yet powerful trick to deal with such PDEs (see also \cite{RasulovRaimovaMascagni2010}). For simplicity we will consider only the one dimensional case. The first trick, is to rescale the problem, by assuming the payoff is bounded. Then we may consider the following function $v:= u/||\psi||_{\infty}$, such a function satisfies the PDE 
\begin{align}
	\label{eq:CVA-PDE-aux}	
	\frac{\partial v}{\partial t}+\cL v + r_{0}v+r_{1}v^{+}=0 
\end{align}
with terminal value bounded above by $1$. Let us consider instead a PDE of the form
\begin{align*}
	\dfrac{\partial v}{\partial t}+\cL v + c\left(F(v) -v \right)=0,
	\qquad 
	v(x,T)=\psi/\|\psi\|_{\infty}
\end{align*}
where $F$ is a polynomial, hence this equation is of the type considered previously. The great insight in \cite{Labordere2012} was to use a polynomial to approximate the semilinear component $v^{+}$. Therefore we have transformed the PDE from one outside the capabilities of marked branching diffusions, into one which can be. \cite{Labordere2012} uses a polynomial of order $4$ to do this, which provides a good approximation as can be seen in Fig.~\ref{Fig:CVAPlot}.
Thus the PDE we want to calculate is,
\begin{align*}
\dfrac{\partial v}{\partial t}+\cL v + c\left(0.0586+0.5v +0.8199v^{2}-0.4095v^{4} -v \right)=0,
 ~ \quad 
v(x,T)=\psi/\|\psi\|_{\infty} \, .
\end{align*}
Since we have constant coefficients and a bounded terminal condition, this PDE easily satisfies Assumption \ref{Assump: Marked Branching Diffusion}, hence we can use the stochastic representation. Remark the negative coefficient associated to the 4th-power term; the classical Assumption \ref{assump:ClassicalKPP} does not hold here and hence the motivation behind Marked Branching Diffusions.

Of course, there does not seem much scope for PDD here since the domain is rather small. However, in higher dimensional problems (basket options), the PDD method can very easily provide the computational advantages highlighted in earlier sections.

\subsection{Further Generalisations: Aged marked branching diffusions}
\label{sec:Agemarked}

In \cite{LabordereEtAl2016}, the authors generalize the representation to a wider class of semilinear PDEs through so called \emph{age-marked branching diffusions}, still building on the paradigm of semilinear PDEs with polynomial source terms $f$. Consider for some integer $m \ge 0$ a set $L \subset \bN^{m+1}$ and a sequence of functions $(c_{l})_{l \in L}$ and $(h_{i})_{i=1, \cdots, m}$, where $c_{l}: [0,T]\times \bR^{d} \rightarrow \bR$ and $h_{i}: [0,T]\times \bR^{d} \rightarrow \bR$. Let the source term function $f: [0,T]\times \bR^{d}\times \bR\times \bR^{d} \rightarrow \bR$ be of the form 
\begin{align*}
f(t,\bfx,y,z)=\sum_{l=(l_{0}, \dots, l_{m}) \in L} c_{l}(t,\bfx)y^{l_{0}}\prod_{i=1}^{m}(h_{i}(t,\bfx) \cdot z)^{l_{i}} \, .
\end{align*}
With such a function in mind (under some regularity assumptions) \cite{LabordereEtAl2016} provides a stochastic representation for PDEs of the form,
\begin{align*}
\frac{\partial u}{\partial t}+\cL u + f(\cdot,u,\nabla u)=0 \, , ~ ~ \text{with }~ u(\cdot,T)=\psi(\cdot),
\quad t\in[0,T], ~ \bfx \in \bR^{d} \, , 
\end{align*}
where $\psi: \bR^{d} \rightarrow \bR$ is a bounded Lipschitz function. The assumptions required for a stochastic representation are more technical than the previous case considered. We thus omit them here and point the reader to \cite{LabordereEtAl2016} for the details. Under this set up we can handle systems of PDE and also popular PDEs such as Burger's equation. In fact, using the polynomial approximation trick discussed in Section \ref{sec:obstaclePDE}, we can handle PDEs containing $f(\cdot,u, \nabla u)= u|\nabla u|^2$ terms, which appear in many applications for instance, to construct harmonic homotopic maps \cite{Struwe1996}; in the theory of ferromagnetic materials through the Landau-Lifschitz-Gilbert equation or models of magnetostriction.

\begin{remark}[Fully nonlinear case]
	Recently, \cite{Warin2017} constructed an algorithm from techniques developed in \cite{LabordereEtAl2016} to approximate a fully non-linear PDEs. This work is still experimental (at the present time) and thus we only reference it here for interested readers.
\end{remark}

\subsubsection{Some other PDE examples in data science solvable by branching diffusions}

We focus on branching diffusions, since they are amenable to PDD.
In mathematical biology one often wants to consider the affect on a population of a species given the presence of another species and the interaction between them. Such models have been considered in works such as \cite{DancerDu1994}, see \cite{GobetLiuZubelli2016} as well. In the paper by Escher and Matioc \cite{EscherMatioc2010} they present a model to describe the growth of tumours. The model is a nonlinear two dimensional system with two decoupled Dirichlet problems.

As highlighted in Section \ref{sec:obstaclePDE}, the methodology described in this work can also be used to deal with PDE obstacle problems. Apart from the well-known problem of pricing American type options we point out to the so-called \emph{Travel agency problem} \cite[Section 5.2]{GobetLiuZubelli2016} where a travel agency needs to decide when to offer travels depending on currency and weather forecast.

If one returns to finance, the approximation method highlighted in Section \ref{sec:obstaclePDE} combined with the age-marked branching diffusions can be used to tackle the problems in actuarial science, namely the ``Reinsurance and Variable Annuity'' problem (see \cite{GuyonLabordere2013}) or the pricing of contracts with different borrowing/lending rates \cite[Section 3.3]{CrisanManolarakis2010}.


\section{An Invitation to More Challenging Problems}
\label{Invitation}

In this Section, we shortly review probabilistic representations of nonlinear problems based on Forward-Backward SDE (FBSDEs) for which branching arises as a particular case. FBSDEs provide stochastic representations for a wide class of PDEs, wider than branching techniques allow, including systems of (fully nonlinear, parabolic or elliptic) PDEs (or BVPs) combining any type of mixed boundary condition. 

Consider the system of $K\geq 1$ quasilinear coupled PDEs in ${\mathbb R}^d \times [0,T]$, $d\geq 1$:
\begin{equation}
\label{F:Nonlin_Sys}
\left\{
\begin{array}{rl}
\frac{\partial u^{(1)}}{\partial t}({\bf x},t) + {\cal L}u^{(1)}({\bf x},t) +
f^{(1)}(t,{\bf x},u^{(1)},\ldots,u^{(K)},\nabla u^{(1)},\ldots,\nabla u^{(K)})& =0, 
\\
u^{(1)}({\bf x},T)&= \psi^{(1)}({\bf x})\\
&\vdots \\
\frac{\partial u^{(K)}}{\partial t}({\bf x},t) + {\cal L}u^{(K)}({\bf x},t) +
f^{(K)}(t,{\bf x},u^{(1)},\ldots,u^{(K)},\nabla u^{(1)},\ldots,\nabla u^{(K)})&=0, 
\\
u^{(K)}({\bf x}, T)&= \psi^{(K)}({\bf x}),
\end{array}
\right. .
\end{equation}
The associated SDE to $\cL$, ${\bf X}_t$, is as in equation \eqref{F:Milstein_sys} (with the same requirements on the coefficients as there). Under the proper conditions on the coefficients in (\ref{F:Nonlin_Sys}) and adequate smoothness of the $u^{(1)}({\bf x},t),\ldots,u^{(K)}({\bf x},t)$  the pointwise solution of the system at $({\bf x}_0,0)$ is given by $(1\leq k\leq K)$
\begin{align}
\label{F:BSDE_Expectation}
u^{(k)}({\bf x}_0,0)
&= {\mathbb E}\Big[ \psi^{(k)}({\bf X}_T) 
\\
\nonumber
&+ \int_0^T f^{(k)}\big( s,{\bf X}_s,Y^{(1)}_s,\ldots,Y^{(K)}_s,{\bf Z}_s^{(1)}\sigma^{-1}({\bf X}_s , s),\ldots,{\bf Z}_s^{(K)}\sigma^{-1}({\bf X}_s , s)\big)ds\Big],
\end{align}
where the processes $Y^{(1)}_t,\ldots,Y^{(K)}_t,{\bf Z}_t^{(1)},\ldots,{\bf Z}_t^{(K)}$ are the solution to the the system of FBSDEs (see \cite[chapter VII]{gobet2016monte})
\begin{equation}
\label{F:BSDEs}
\left\{
\begin{array}{rl}
Y_t^{(1)}&= \psi^{(1)}({\bf X}_T) -\int_t^T {\bf Z}^{(1)}_s\cdot d{\bf W}_s
\\
&\quad + \int_t^T f^{(1)}\big( s,{\bf X}_s,Y^{(1)}_s,\ldots,Y^{(K)}_s,{\bf Z}_s^{(1)}\sigma^{-1}({\bf X}_s , s),\ldots,{\bf Z}_s^{(K)}\sigma^{-1}({\bf X}_s , s)\big)ds ,\\
&\vdots\\
Y_t^{(K)}&= \psi^{(K)}({\bf X}_T) -\int_t^T {\bf Z}^{(K)}_s\cdot d{\bf W}_s
\\
&
\quad + \int_t^T f^{(K)}\big( s,{\bf X}_s,Y^{(1)}_s,\ldots,Y^{(K)}_s,{\bf Z}_s^{(1)}\sigma^{-1}({\bf X}_s , s),\ldots,{\bf Z}_s^{(K)}\sigma^{-1}({\bf X}_s , s)\big)ds .
\end{array}
\right.
\end{equation}
One can show via It\^o's formula that $Y_t^{(k)}=u^{(k)}({\bf X}_t , t)$ and ${\bf Z}^{(k)}_t=\nabla u^{(k)}({\bf X}_t , t)\sigma({\bf X}_t , t)$, $(1\leq k\leq K)$ --- this lifts the apparent inconsistency of having more processes than equations in (\ref{F:BSDEs}). As with linear parabolic equations, the terminal condition can be transformed into an initial one by reversing time, i.e. by letting ${\hat t}\shortrightarrow T-t$ and $\partial/\partial {\hat t}\shortrightarrow -\partial/\partial t$. System (\ref{F:Nonlin_Sys}) accommodates many important applications in biochemistry, stochastic control, finance and physics \cite[chapter VII]{gobet2016monte}. Multiple extensions of the above system are possible, including different generators for each of the equations in the system (i.e. different ${\cal L}^{(1)},\ldots,{\cal L}^{(K)}$), as well as BCs of various types. The complexity of the probabilistic representation grows accordingly.

In order to simplify the discussion and focus on difficulties, let $K=1$ in the remainder of this Section (thus $Y_t=Y^{(1)}_t$, $\psi=\psi^{(1)}$, $f=f^{(1)}$), and $f=f(t,{\bf x},u)$. It can be shown 
\begin{equation}
\label{F:Conditional_expectation}
Y_t=u({\bf X}_t , t)= {\mathbb E}\Big[ \psi({\bf X}_T) + \int_t^T f(s,{\bf X}_s,Y_s)ds \,\big|\, {\bf X}_t\Big].
\end{equation}
The manipulation of the conditional expectation (\ref{F:Conditional_expectation}) requires the introduction of filtered probability spaces, which lies beyond the scope of this paper. The numerical treatment is equally delicate. After a canonical time discretization over a uniform time partition with step $h=T/N$, we have
\begin{align}
\label{eq:BSDE-discretization}
Y_{t_N}=\psi({\bf X}_T),\qquad \qquad Y_{t_i}&= {\mathbb E}\Big[ Y_{t_{i+1}} + h f(t_i,{\bf X}_{t_i},Y_{t_{i+1}}) \,\big|\, {\bf X}_{t_i}\Big]
\\
\nonumber
\Leftrightarrow  u({\bf X}_{t_i},t_i)&= {\mathbb E}\Big[ Y_{t_{i+1}} + h f\big(t_i,{\bf X}_{t_i},u({\bf X}_{{t_{i+1}}},{t_{i+1}}) \big) \,\big|\, {\bf X}_{t_i}\Big].
\end{align}
Here we emphasize that to compute $u({\bf x},0)=Y_0$ one needs, at each time step $t_i$, to compute an approximation of $u({\bf x}, t_i)$ for $x$ over the domain of the PDE. Comparing with the result when branching is possible, see equation \eqref{Eq:Classic Branching with Boundary}, the gain of branching over the general FBSDE machinery is clear. Nonetheless, branching techniques have their limits and the general case must, in principle, be tackled through FBSDEs and its variants. How to carry out a PDD methodology within the framework of FBSDEs that still retains PDD's characteristic scalability is an open question.

In general, the timestep-per-timestep backward iteration and the computation and/or approximation of the conditional expectation in \eqref{eq:BSDE-discretization} can be done in many ways, for instance, by resorting to Bellman's dynamic programming principle, whereby the solution is projected on a functional basis by regression (see \cite[chapter VIII]{gobet2016monte}). Overviews on numerics of FBSDE can be found in \cite{BouchardElieTouzi2009}, \cite{CrisanManolarakis2010} and PDE inspired problems dealt by FBSDEs can be found in \cite{FreiDosReis2013}, \cite{LionnetdosReisSzpruch2015}, \cite{LionnetdosReisSzpruch2016}.

\subsubsection*{A brief remark on FBSDEs and stochastic representations}

We close this section with comments on various probabilistic representations not mentioned elsewhere in this paper. 

There is theoretical work relating FBSDEs to the Navier-Stokes equations \cite{cruzeiro2009navier}. While, in general, transport equations do not have a stochastic representation, there are some important exceptions, of which we cite two: the one-dimensional telegraph equation \cite{AcebronRibeiro2016}; and the the Vlasov-Poisson system of equations. For the latter, it was proposed in \cite{Mendes2010} a probabilistic representation in the Fourier space.
 
In \cite{BossyEtAl2015} the authors deal with representations (and Monte Carlo simulation) for nonlinear divergence-form elliptic Poisson-Boltzmann PDE over the whole $\bR^3$. Quasilinear parabolic PDEs (multidimensional and systems) admit a stochastic representation in terms of FBSDEs \cite{Peng1991}, \cite{PardouxPeng1992} and many extensions exist ranging from representations for obstacle PDE  problems \cite{ElKarouiKapoudjianPardouxEtAl1997} to representations of fully nonlinear elliptic and parabolic PDEs \cite{CheriditoSonerTouziEtAl2007}. See as well \cite{FreiDosReis2013}, \cite{LionnetdosReisSzpruch2015} for stochastic representations for reaction-diffusion PDEs, PDEs with quadratic gradients terms and of Burger's type nonlinearities. Complex valued PDEs and the FBSDE connection are handled in \cite{xu2015complex}.

Systems of fully nonlinear systems of second-order PDEs (i.e. including possible nonlinearities on the highest derivatives i.e. $f$ can depend on ``$\Delta u$'' terms) have recently found a representation in terms so-called 2BSDEs \cite{CheriditoSonerTouziEtAl2007}; numerical schemes have already been proposed \cite{fahim2011probabilistic}. Lastly, much of the theory (on stochastic representations) can be extended to frameworks where the coefficients $b$ and $\sigma$ depend on $u$ and $\nabla u$ (see \cite{MaYong1999}).

\begin{acknowledgements}
F.~Bernal acknowledges funding from \emph{Centre de Math\'ematiques Appliqu\'ees} (CMAP), \'Ecole Polytechnique.

G.~dos Reis gratefully thanks the partial support by the \emph{Funda\c c\~ao para a Ci\^encia e a Tecnologia} (Portuguese Foundation for Science and Technology) through the project UID/MAT/00297/2013 (\emph{Centro de Matem\'atica e Aplica\c c\~oes}).

G.~Smith was supported by The Maxwell Institute Graduate School in Analysis and its
	Applications, a Centre for Doctoral Training funded by the UK Engineering and Physical
	Sciences Research Council (grant [EP/L016508/01]), the Scottish Funding Council, Heriot-Watt
	University and the University of Edinburgh.
\end{acknowledgements}

\end{document}